%% file: main.tex
\documentclass[]{llncs}
\usepackage{amsmath}
\usepackage{amssymb}
\usepackage{txfonts}
\usepackage{stmaryrd} 
\usepackage{rrproof}
\usepackage{paralist}

\include{macros}

\shortfalse

\usepackage{float}    
\floatstyle{boxed}
\restylefloat{figure}
\floatstyle{ruled}
\restylefloat{table}
\newcommand{\IntGeo}{\ensuremath{\sn{Int^*}/\mt{Geo}}}
\renewcommand{\defterm}[1]{\emph{#1}}
\newcommand{\RLS}[1]{{#1}_{\rel}}  
\newcommand{\TU}[1]{{#1}^{\bullet}}
\newcommand{\ol}[1]{\bar{#1}}
\newcommand{\ox}{\ol{x}}
\newcommand{\oy}{\ol{y}}

\spnewtheorem{conj}{Conjecture}{\bf}{\it}
\spnewtheorem*{theorem*}{Theorem}{\bf}{\it}

\title{Translating Labels to Hypersequents for Intermediate Logics with Geometric Theories}
\author{Robert Rothenberg\inst{1}\inst{2}
\institute{School of Computer Science, University of St Andrews, Fife KY16 9SX, Scotland
   \email{rr@cs.st-andrews.ac.uk} \\
\and
Interactive Information, Ltd., 3 Lauriston Gdns, Edinburgh EH3 9HH, Scotland
 }
}

\begin{document}

\maketitle

\begin{abstract}
We give a procedure for translating geometric Kripke frame axioms
into structural hypersequent rules for the corresponding
intermediate logics with geometric theories (\IntGeo) that admit
weakening, contraction and cut. We give a procedure
for translating labelled sequents in the corresponding logic to
hypersequents that share the same linear models (which correspond
to G\"odel-Dummett logic). We prove that labelled proofs of formulae
for logics in \IntGeo~ can be translated into hypersequent proofs
that in special cases may use the linearity rule, which corresponds to the
well-known communication rule for G\"odel-Dummett logic.
\end{abstract}

\section{Introduction}
Since the introduction of the sequent calculus by Gentzen \cite{Gent1935},
the syntactic elements of the calculus have been extended to give
calculi for various non-classical logics.  The formal relationships between two
common extensions, labelled sequents and hypersequents, have not been
examined in the literature, besides translations between \emph{specific}
calculi. (These will be discussed in \S\ref{sec:related} below.)

We give a method for translating labelled sequent proofs into
hypersequent proofs, for a subset of intermediate logics with
geometric frame axioms. The resulting hypersequent proofs are
sometimes in stronger logics based on G\"odel-Dummett Logic (\sn{GD} \cite{Dumm1959-JSL})
rather than Intuitionistic Logic (\sn{Int}).

Labelled sequent calculi, apparently introduced in \cite{Kang1957},
contain formulae that are annotated with labels, and often the
sequents themselves are annotated with terms that indicate the
relationships between labels.  Hypersequents, generally attributed to \cite{Avro1987} (though they have occurred earlier, e.g. \cite{Beth1959} and \cite{Pott1983}),
are lists or multisets of sequents.

Developing a formal translation between proof systems is a topic of
interest.  The obvious reasons for doing so are to allow one to
separate interface from implementation in automated proof assistants
(especially where one formalism is more conducive to automation), and
to translate proofs of meta properties, such as cut-elimination or
interpolation, into alternative formalisms.

A less obvious reason for developing translations are to gain a better
understanding of the \emph{meaning} of particular syntactic features
that proof systems extend with respect to sequent calculi, where the
``meaning'' of a syntactic feature is the function that it plays
in the inference rule.

This is useful for developing new notations which can combine multiple
syntactic features.  Such a notation can be used to develop new
extensions to sequent calculi, or to develop a formal hierarchy of the
relative strength of proof systems.

Labelled calculi can be seen as an alternative notation for other
formalisms, where the locations of formulae in a structure are encoded
as labels, and the relationships between these locations
are encoded as relational formulae.  This location
information about the data structure as part of the object language of
a labelled calculus.  For example, labels can name the components
of the hypersequent-like structure,
$$
\underbrace{ A, \ldots}_{x} ~|~ \underbrace{A, B, \ldots}_{y} ~|~ \underbrace{A ,B, C \ldots}_{z}
$$
that can then be translated into a kind of labelled sequent:
$\lb{x,y,z,\ldots}{A},\lb{y,z\ldots}{B},\lb{z,\ldots}{C}$.
Relations can be added to encode relationships (such as subset relations)
between components:
$x\rel y,y\rel z ; \lb{x,\ldots}{A},\lb{y,\ldots}{B},\lb{z,\ldots}{C}$.
Labels and relations can be used to reduce the complexity of
data structures, e.g. the above structure may be easier to search for
formulae in than the original structure.
Likewise, a decoding of relationships may also be useful, e.g.
to translate a structure into something that may be easier to search
using parallel algorithms.

\subsection{Related Work}
\label{sec:related}

The relationship between labelled sequents and hypersequents has been
a folkloric one in proof theory, with no published formal comparisons
that the author is are aware of, beyond those for specific calculi.
Much of the work has been for systems based on the modal logic \sn{S5}
\cite{Avro1996}, \cite{Rest2006-Ms} and \cite{GiraP2005-IMLA} (the
latter work also connects systems for the Logic of Strong Negation,
\sn{N3} \cite{Nels1973-SL}).  Slightly more general work on systems
of modal logic can be found in \cite{Melvin200785}.
Work connecting specific hypersequent
and labelled calculi for Abelian Logic (\sn{A} \cite{MeyeSlan1989})
and \Luka Logic (\textbf{\L} \cite{Luka1930a})
is given in \cite{MetcOlivGabb2005}.

General work on deriving a relational semantics, which can be used as the
basis for labelled calculi, from Hilbert- or Gentzen-style calculi
(which presumably can be extended to hypersequents) is given in
\cite{Gabb96-Book}.

Work on using the form of Hilbert axioms to determine the
kind of calculus required (e.g. sequent or hypersequent) for a
cut-free proof system, and on obtaining structural rules corresponding
from those axioms, is given in \cite{CiabGalaTeru2008-IEEESLCS}.

Work on translating some Kripke frame axioms (what we call ``geometric
theories'' in this paper) into structural rules that admit weakening,
contraction and cut for a \sn{G3}-style labelled calculus has been
done by various authors, notably \cite{Simp1994-PhD}, \cite{Viga2000}
and \cite{Negr2005-ASL}.  The latter work was adapted to general work
on translating between hypersequents and labelled sequents for logics
in \IntGeo~ in \cite{Roth2010-PhD}, and is used as a basis for parts
of this paper---in particular, a method was given for translating
labelled sequent proofs for logics in \IntGeo~into simply labelled
proofs (an alternative notation for hypersequent proofs discussed later)
for a corresponding calculus augmented with a form of the
communication rule from \cite{Avro1991,Avro1996}.

Discussions of syntactic extensions to sequent calculus can be found in
\cite{Paol2002-Book}. 

\section{Preliminaries}

We give a brief overview of the class of logics, \IntGeo, along with
labelled sequents, hypersequents, and simply labelled sequents, which
will be used to give calculi for logics in that class.

\subsection{General Notation}

Uppercase Roman letters such as $A,B,C$ will denote arbitrary logical
formulae.  The letters $P,Q,R$ will denote atomic formulae.  Uppercase
Greek letters such as $\Gm,\Dl,\Pi,\Phi,\Psi$ will denote
(possibly empty) multisets of arbitrary logical formulae.

When introducing a rules for a calculus, double lines will be used to
indicate that a rule is invertible. (However, a single line does
not necessarily mean that a rule is not invertible.)

In proofs or proof fragments, an exponent of $n$ on a rule name
indicates $n$ applications of the rule, an exponent of $*$ indicates
$0$ or more applications of the rule, and an exponent of $+$ indicates
$1$ or more applications of the rule.

\subsection{Intermediate Logics with Geometric Theories}

Intermediate Logics (\sn{Int^*}) are (propositional) logics between
Intuitionistic Logic (\sn{Int}) and Classical Logic (\sn{Class}) that
can be obtained by adding additional axioms to \sn{Int}.
(See Table~\ref{tbl:char} for well-known examples.)  Below, we give
a semantic characterisation of a subclass of them, \IntGeo, that we
call \defterm{Intermediate Logics with Geometric Theories}.

\begin{definition}
\label{definition:kripke_semantics_for_int}
An \defterm{Intuitionistic Kripke Frame}
is a structure $\struct{\set{W},R}$ where
$W$ is a set of atomic \defterm{points}, $R$ is preordered
binary relation on $W$. An \defterm{Intuitionistic Krip\-ke Model}
$\mf{M}$ is an Intuitionistic Krip\-ke Frame extended with $D$, a function
from points to sets of atomic formulae, which is {monotonic} w.r.t.
$R$---i.e., for all $x,y\in W$, if
$Rxy$, then $D(x)\subseteq D(y)$.
A \defterm{forcing relation} $\mf{M},x\force A$ for propositional formulae
is defined as follows:
\begin{enumerate}
\item $\mf{M},x\force P$ iff $P\in D(x)$ for all $x\in W$.

\item $\mf{M},x\nforce\bot$, i.e. $\bot\notin D(x)$ for all $x\in W$.\label{item:bot}

\item $\mf{M},x\force A\land B$ iff $\mf{M},x\force A$ and $\mf{M},x\force B$.

\item $\mf{M},x\force A\lor B$ iff either $\mf{M},x\force A$ or $\mf{M},x\force B$.

\item $\mf{M},x\force A\imp B$ iff for all $y$ such that $Rxy$,
$\mf{M},y\force A$ implies $\mf{M},y\force B$.

\ifshort
\else

\item $\mf{M},x\force\neg A$ iff for all $y$ such that $Rxy$,
$\mf{M},y\nforce A$.

\item $\mf{M},x\force\top$ for all $x\in W$.

\fi
\end{enumerate}
If $\mf{M},x\force A$ for all $x\in W$, then we write simply that
$\mf{M}\vDash A$.
When $\mf{M}$ is obvious from the context, we write $x\force A$.
\end{definition}

Models for many logics in \sn{Int^*} can be obtained by extending the
frames of an Intuitionistic Kripke Model $\mf{M} =
\struct{\set{W},R,D}$ with additional axioms on $R$.  For many
well-known logics, such as those in Table~\ref{tbl:char}, the frame
axioms are \defterm{geometric implications}---that is, they are of the
form $\forall\ox.(A\imp B)$, with restrictions given in
Definition~\ref{def:geometric}.
The logics that correspond to such models are said to be
in the class \IntGeo.  These logics are of interest because the
structural rules that correspond to their characteristic frame axioms
can be added to \sn{G3}-style labelled sequent calculi without
affecting the admissibility of the standard structural rules
\cite{Negr2005-ASL}, as will be discussed below.

\begin{definition}[Geometric Formulae and Implications]\label{def:geometric}
A \defterm{geometric formula} is defined inductively:
\begin{enumerate}
\item An atomic formula $P$ is a geometric formula.
\item $\bot$ is a geometric formula.
\item $\top$ (that is, $\bot\imp\bot$) is a geometric formula.
\item $A\lor B$ is a geometric formula iff $A$ and $B$ are geometric formulae.
\item $A\land B$ is a geometric formula iff $A$ and $B$ are geometric formulae.
\item $\exists\bar{x}.A$ is a geometric formula iff $A$ is a geometric formula.
\end{enumerate}
A \defterm{geometric implication} is a formulae of the form
$\forall\bar{x}.(A\imp B)$, where $A$ and $B$ are geometric formulae.
\end{definition}

\begin{remark} Reflexivity and transitivity axioms are geometric implications.
\end{remark}

\begin{table}[h]
\caption{Some well-known logics with their characteristic axioms and frame axioms.}\label{tbl:char}
\centering
\begin{tabular}{c c c}
 \textbf{Logic} &
 \textbf{Axiom} &
 \textbf{Frame Axiom}
\\
\\
 Jankov-De Morgan (\sn{Jan}) &
 $\neg A\lor\neg\neg A$  &
 $\forall xy\in W ~.~ \exists z\in W ~.~ Rxz \lor Ryz$
\\
\\
 G\"odel-Dummett (\sn{GD})~~  &
 $(A\imp B)\lor(B\imp A)$ &
 $\forall xy\in W ~.~ Rxy\lor Ryx$
\\
\\
 Bounded-Depth of 2 (\sn{BD_2})~~ &
 $B\lor(B\imp(A\lor\neg A))$ &
 ~~ $\forall xyz\in W~.~ Rxy\land Ryz\ \imp Ryx\lor Rzy$
\\
\\
 Classical (\sn{Class})  &
 $A\lor \neg A$ &
 $\forall xy\in W ~.~ Rxy\imp Ryx$
\end{tabular}
\end{table}

\begin{lemma}[Pointed Models \cite{Mint2000} \S 7.2]\label{lem:rooted}
Let $\mf{M} = \struct{W,R,D}$ and $\mf{M}' = \struct{W',R',D'}$ be
Kripke models for a logic in \sn{Int^*}, such that
$\mf{M}\vDash A$ iff $\mf{M}'\vDash A$.
Then if $\mf{M}'$ is \defterm{pointed}, i.e.
$\exists x\in W'$ (a \defterm{distinguished point}) s.t.
$\forall y\in W'$, $Rxy'$, then $R'$ is a partial order,
i.e. it is also anti-symmetric.
\end{lemma}

\begin{remark} In \cite{DyckNegr2011}, the frame axiom for
\sn{Jan} (Table~\ref{tbl:char}) is given for pointed models,
i.e. $\forall wxy\in W . Rwx\land Rwy\imp\exists z\in W . Rxz \lor Ryx$.
However, both versions are interderivable by Lemma~\ref{lem:rooted}.
\end{remark}

\subsection{Labelled Sequent Calculi}

Labelled sequents are an extension of Gentzen-style sequents, where
the logical formulae are annotated with (atomic) labels,
e.g. $\lb{x}{(A\lor B)}$.  It is common in contemporary systems, such
as \cite{Negr2005-ASL}, that the sequents are also annotated with a
collection of (binary) relations, called \defterm{relational formulae},
between labels, e.g. $x\rel y$. (Such systems are generally used to
reason about a logic's corresponding relational models.)  More
expressive labelled calculi are discussed in \cite{Gabb96-Book,Blac2000-LD}.

We denote labelled sequents as $\Sg;\nlb{\Gm}\seq\nlb{\Dl}$,
where $\Sg$ is an arbitrary multiset of relational formulae, and
underlined multiset variables are multisets of formulae with arbitrary labels.
(Labelled multiset variables, e.g. $\lb{x}{\Gm}$, denote multisets
of formulae with the same label.)
$x\#\nlb{\Gm}$ denotes that the label $x$ is fresh for
$\nlb{\Gm}$.
Multiple occurrences of a labelled formula may
 be abbreviated by concatenating the labels, e.g.
$\lb{xxy}{A} = \lb{x}{A},\lb{x}{A},\lb{y}{A}$. Multisets
of labelled formulae may similarly be abbreviated.

Symmetric relational formulae such as $x\rel y,y\rel x$
may be abbreviated as $x=y$, and transitive pairs
of relational formulae such as $x\rel y,y\rel z$ may be abbreviated
as $x\rel y\rel z$.

The semantics for labelled sequents is given in
Definition~\ref{def:rls_semantics}.  The vocabulary describing sequents from \cite{TroeASchwH2000} is extended naturally
for labelled sequents. A calculus \sn{G3I}
\cite{Negr2005-ASL,DyckNegr2011} for \sn{Int} is given in
Fig.~\ref{figure:G3I}.

\begin{definition}[Semantics of Labelled Sequents]\label{def:rls_semantics}
Let $\mf{M} = \struct{W,R,D}$ be a Kripke model for a logic in
\IntGeo.  Then $\mf{M}\vDash \Sg;\nlb{\Gm}\seq\nlb{\Dl}$ iff
for each $w\in\mt{lab}(\Sg,\nlb{\Gm},\nlb{\Dl})$, there exists a
(not necessarily unique) $\hat{w}\in W$, such that
the \defterm{consistency of $\Sg$ with $R$}---i.e., for all $x\rel y\in\Sg$, $R\hat{x}\hat{y}$---implies
either $\mf{M}\nvDash\mland\nlb{\Gm}$ or $\mf{M}\vDash\mlor\nlb{\Dl}$,
where $\mf{M}\vDash\lb{x}{A}$ iff $\mf{M},\hat{x}\force A$.
\end{definition}

\begin{figure}[h]
\centering
\begin{tabular}{c c}
$$
\infer[\RLS{Ax}]{x\rel y,\Sg;\lb{x}{P},\nlb{\Gm}\seq\lb{y}{P},\nlb{\Dl}}{
}
$$
&
$$
\infer[{L\bot}]{\Sg;\lb{x}{\bot},\nlb{\Gm}\seq\nlb{\Dl}}{
}
$$
\\ \\
\ifshort
\else
$$
\infer[{L\land}]{\Sg;\nlb{\Gm},\lb{x}{(A\land B)}\seq\nlb{\Dl}}{
 \Sg;\nlb{\Gm},\lb{x}{A},\lb{x}{B}\seq\nlb{\Dl}
}
$$
&
$$
\infer[{R\land}]{\Sg;\nlb{\Gm}\seq \lb{x}{(A\land B)},\nlb{\Dl}}{
 \Sg;\nlb{\Gm}\seq \lb{x}{A},\nlb{\Dl} &
 \Sg;\nlb{\Gm}\seq \lb{x}{B},\nlb{\Dl}
}
$$
\\ \\
\fi
$$
\infer[{L\lor}]{\Sg;\nlb{\Gm},\lb{x}{(A\lor B)}\seq\nlb{\Dl}}{
 \Sg;\nlb{\Gm},\lb{x}{A}\seq\nlb{\Dl} &
 \Sg;\nlb{\Gm},\lb{x}{B}\seq\nlb{\Dl}
}
$$
&
$$
\infer[{R\lor}]{\Sg;\nlb{\Gm}\seq \lb{x}{(A\lor B)},\nlb{\Dl}}{
 \Sg;\nlb{\Gm}\seq \lb{x}{A},\lb{x}{B},\nlb{\Dl}
}
$$
\\ \\
\multicolumn{2}{c}{$$
\infer[\RLS{L\imp}]{x\rel y,\Sg;\lb{x}{(A\imp B)},\nlb{\Gm}\seq\nlb{\Dl}}{
 x\rel y,\Sg;\lb{x}{(A\imp B)},\nlb{\Gm}\seq\nlb{\Dl},\lb{y}{A} &
 x\rel y,\Sg;\lb{x}{(A\imp B)},\lb{y}{B},\nlb{\Gm}\seq\nlb{\Dl}
}
$$}
\\ \\
 &
$$
\infer[\RLS{R\imp}]{\Sg;\nlb{\Gm}\seq\nlb{\Dl},\lb{x}{(A\imp B)}}{
 x\rel y,\Sg;\lb{y}{A},\nlb{\Gm}\seq\nlb{\Dl},\lb{y}{B}
}
$$
\\ \\
&
{where $y\#\Sg,\nlb{\Gm},\nlb{\Dl}$ in \rn{\RLS{R\imp}}}
\\ \\
$$
\infer[{refl}]{\Sg;\nlb{S}}{
  x\rel x,\Sg;\nlb{S}
}$$
&
$$
\infer[{trans}]{x\rel y\rel z,\Sg;\nlb{S}}{
  x\rel z,x\rel y\rel z,\Sg;\nlb{S}
}
$$
\end{tabular}
\caption{The labelled calculus \sn{G3I}. }
\label{figure:G3I}
\end{figure}

\begin{proposition}[\cite{DyckNegr2011}] Labelled sequents of the form
$x\rel y,\Sg;\lb{x}{A},\nlb{\Gm}\seq\lb{y}{A},\nlb{\Dl}$ are derivable
in \sn{G3I}.\label{prop:general_ax_G3I}
\end{proposition}

\begin{lemma}[\cite{Negr2005-ASL,DyckNegr2011}] The weakening, contraction and \rn{cut} rules\label{prop:G3I_struc}

\begin{tabular}{c c c}
\\
$$
\infer[\RLS{LW}]{x\rel y,\Sg;\nlb{\Gm}\seq\nlb{\Dl}}{
 \Sg;\nlb{\Gm}\seq\nlb{\Dl}
}
$$
&
$$
\infer[LW]{\Sg;\lb{x}{A},\nlb{\Gm}\seq\nlb{\Dl}}{
 \Sg;\nlb{\Gm}\seq\nlb{\Dl}
}
$$
&
$$
\infer[RW]{\Sg;\nlb{\Gm}\seq\nlb{\Dl},\lb{x}{A}}{
 \Sg;\nlb{\Gm}\seq\nlb{\Dl}
}
$$
\\
\\
$$
\infer[\RLS{LC}]{x\rel y,\Sg;\nlb{\Gm}\seq\nlb{\Dl}}{
 x\rel y,x\rel y,\Sg;\nlb{\Gm}\seq\nlb{\Dl}
}
$$
&
$$
\infer[LC]{\Sg;\lb{x}{A},\nlb{\Gm}\seq\nlb{\Dl}}{
 \Sg;\lb{x}{A},\lb{x}{A},\nlb{\Gm}\seq\nlb{\Dl}
}
$$
&
$$
\infer[RC]{\Sg;\nlb{\Gm}\seq\nlb{\Dl},\lb{x}{A}}{
 \Sg;\nlb{\Gm}\seq\nlb{\Dl},\lb{x}{A},\lb{x}{A}
}
$$
\\
\end{tabular}
$$
\infer[cut]{\Sg;\nlb{\Gm}\seq\nlb{\Dl}}{
 \Sg;\nlb{\Gm}\seq\nlb{\Dl},\lb{x}{A} &
 \Sg;\lb{x}{A},\nlb{\Gm}\seq\nlb{\Dl}
}
$$
are admissible in \sn{G3I^*}.
\end{lemma}

\begin{remark}[Notation]\label{rem:short_struc} For brevity, multiple instances of \rn{LW\rel}, \rn{LW} or \rn{RW} will be shown in simply as \rn{W}. Similarly, multiple instances of \rn{LC\rel}, \rn{LC} or \rn{RC} will be shown simply as \rn{C}.
\end{remark}

\begin{samepage}
\begin{proposition}\label{prop:rel_rules} The rules
$$
\infer[L\rel]{x\rel y,\Sg;\nlb{\Gm},\lb{x}{A}\seq\nlb{\Dl}}{
 x\rel y,\Sg;\nlb{\Gm},\lb{x}{A},\lb{y}{A}\seq\nlb{\Dl}
}
\qquad
\infer[R\rel]{x\rel y,\Sg;\nlb{\Gm}\seq\nlb{\Dl},\lb{y}{A}}{
 x\rel y,\Sg;\nlb{\Gm}\seq\nlb{\Dl},\lb{x}{A},\lb{y}{A}
}
$$
are admissible in \sn{G3I}.
\begin{proof} Using Proposition~\ref{prop:general_ax_G3I} and \rn{cut}.
\ifshort
\else
$$
\infer[cut]{x\rel y,\Sg;\nlb{\Gm},\lb{x}{A}\seq\nlb{\Dl}}{
 x\rel y,\Sg;\nlb{\Gm},\lb{x}{A}\seq\nlb{\Dl},\lb{y}{A}
 &
 x\rel y,\Sg;\nlb{\Gm},\lb{x}{A},\lb{y}{A}\seq\nlb{\Dl}
}
$$

$$
\infer[cut]{x\rel y,\Sg;\nlb{\Gm}\seq\nlb{\Dl},\lb{y}{A}}{
 x\rel y,\Sg;\nlb{\Gm}\seq\nlb{\Dl},\lb{x}{A},\lb{y}{A}
 &
 x\rel y,\Sg;\nlb{\Gm},\lb{x}{A}\seq\nlb{\Dl},\lb{y}{A}
}
$$

\fi
\end{proof}
\end{proposition}
\end{samepage}

\ifshort
\else
\begin{remark}
The \rn{L\rel} and \rn{R\rel} rules are primitive in the
system \sn{L} for \sn{BiInt} \cite{PintUust2009,DBLP:journals/corr/abs-1101-5445}.
\end{remark}
\fi

\subsubsection{Geometric Rules for Intermediate Logics.}

In \cite{DBLP:journals/mlq/Palmgren02}, it was shown that any set of
geometric implications is {constructively equivalent} to a set
consisting of formulae of the form
$\forall\ox.(A_0\imp\exists\oy.(A_1\lor\ldots\lor A_n))$, where each
$A_i$ is a conjunction of atomic formulae, such as relational formulae.
Formulae in that form, such as the frame
axioms from Table~\ref{tbl:char}, can be translated into rules of the form:
$$
\vcenter{
\infer[]{\overline{A_0},\Sg;\nlb{\Gm}\seq\nlb{\Dl}}{
 \overline{A_1},\overline{A_0},\Sg;\nlb{\Gm}\seq\nlb{\Dl} &
 \ldots &
 \overline{A_n},\overline{A_0},\Sg;\nlb{\Gm}\seq\nlb{\Dl}
}}
$$
where (in an abuse of notation) $\overline{A_i}$ is the multiset of relational formulae in $A_i$, and the variables correspond to labels.
A translation method is given in Definition~\ref{def:deriving_geo_rules}:
\begin{definition}\label{def:deriving_geo_rules}
Given a geometric implication of the form $\forall\ox.(A_0\imp\exists\oy.(A_1\lor\ldots\lor A_n))$, the corresponding \defterm{geometric rule}
can be obtained by
straightforward analysis of the sequent
$\forall\ox.(A_0\imp\exists\oy.(A_1\lor\ldots\lor A_n)),\overline{A_0},{\Gm}\seq{\Dl}$
in a \sn{G3}-style sequent calculus for \sn{Int}, such as \sn{G3i}
\cite{TroeASchwH2000} using \begin{inparaenum}[(1)]
\item  \rn{L\forall};
\item \rn{L\imp};
\item \rn{R\land};
\item \rn{L\exists}; and
\item \rn{L\lor}.
\end{inparaenum}
\end{definition}

\begin{figure}
$$
\infer[dir]{\Sg;\nlb{\Gm}\seq\nlb{\Dl}}{
 x\rel{z},y\rel z,\Sg;\nlb{\Gm}\seq\nlb{\Dl}
}
\quad
\infer[{lin}]{\Sg;\nlb{\Gm}\seq\nlb{\Dl}}{
 x\rel y,\Sg;\nlb{\Gm}\seq\nlb{\Dl} & y\rel x,\Sg;\nlb{\Gm}\seq\nlb{\Dl}
}
$$

$$
\infer[bd_2]{x\rel y\rel z,\Sg;\nlb{\Gm}\seq\nlb{\Dl}}{
 x = y\rel z,\Sg;\nlb{\Gm}\seq\nlb{\Dl} &
 x\rel y = z,\Sg;\nlb{\Gm}\seq\nlb{\Dl}
}
\quad
\infer[{sym}]{x\rel y,\Sg;\nlb{\Gm}\seq\nlb{\Dl}}{
 x=y,\Sg;\nlb{\Gm}\seq\nlb{\Dl}
}
$$
\caption{Extension rules to \sn{G3I} some well-known logics. ($z\#\Sg;\nlb{\Gm},\nlb{\Dl}$ in \rn{dir}.)}\label{fig:char_rules}
\end{figure}

In \cite{Negr2003-AML} it was shown that {geometric rules}
can be added to \sn{G3}-style calculi without losing the
admissibility of weakening, contraction and \rn{cut}.
This allowed the development of labelled sequent frameworks for various
non-classical logics in \cite{Negr2005-JPL} and \cite{Negr2005-ASL}.
The corresponding rules to frame axioms from Table~\ref{tbl:char} are
in Fig.~\ref{fig:char_rules}. We denote \sn{G3I} augmented
with arbitrary geometric rules such as those from Table~\ref{tbl:char}
as \sn{G3I^*}.

\subsection{Hypersequent Calculi}

A hypersequent is a non-empty multiset of sequents, called its
\defterm{components}, and is written as
$\Gm_1\seq\Dl_1~|~\ldots~|~\Gm_n\seq\Dl_n$.  Hypersequent rules are
written similarly to sequent rules, with calligraphic Roman letters
$\hG,\hH$ used as metavariables to denote a (possibly empty)
multiset of side components in a hypersequent rule.

\begin{definition}[Semantics of Hypersequents]\label{def:hs_semantics}
Let $\mf{M}$ be a model for a logic in
\IntGeo. Then $\mf{M}\vDash\Gm_1\seq\Dl_1~|~\ldots~|~\Gm_n\seq\Dl_n$ if
there exists $1\le i\le n$ such that $\mf{M}\vDash\Gm_i\seq\Dl_i$, i.e.
either $\mf{M}\nvDash\mland\Gm_i$ or $\mf{M}\vDash\mlor\Dl_i$.
($\Gm_i\seq\Dl_i$ is called the \defterm{distinguished component}.)
\end{definition}

The vocabulary describing sequents from \cite{TroeASchwH2000} is
extended naturally for hypersequents.
Rules of hypersequent calculi can be classified as either
\defterm{internal rules} (rules which have only one active component
in each premiss and one principal component in the conclusion), and
\defterm{external rules}, which are not internal rules.
The \defterm{standard external rules} are \rn{EW} and \rn{EC}
(see Fig.~\ref{figure:H3ip}).

The hypersequent calculus \sn{HG3ipm} \cite{Roth2010-PhD} for \sn{Int}
given in Fig.~\ref{figure:H3ip} was obtained from a multisuccedent
variant of \sn{G3ip} \cite{TroeASchwH2000} by adding side components
to the rules and the standard external rules to the calculus.

\begin{figure}
\centering
\begin{tabular}{c c}
$$
\infer[Ax]{\hH~|~P,\Gm\seq\Dl,P}{
}
$$
&
$$
\infer[L\bot]{\hH~|~\bot,\Gm\seq\Dl}{
}
$$
\\ \\
\ifshort
\else
$$
\infer[L\land]{\hH~|~\Gm,A\land B\seq\Dl}{
  \hH~|~\Gm,A, B\seq\Dl
}
$$
&
$$
\infer[R\land]{\hH~|~\Gm\seq A\land B,\Dl}{
  \hH~|~\Gm\seq A,\Dl &
  \hH~|~\Gm\seq B,\Dl
}
$$
\\ \\
\fi
$$
\infer[L\lor]{\hH~|~\Gm,A\lor B\seq\Dl}{
  \hH~|~\Gm,A\seq \Dl &
  \hH~|~\Gm,B\seq \Dl
}
$$
&
$$
\infer[R\lor]{\hH~|~\Gm\seq A\lor B,\Dl}{
  \hH~|~\Gm\seq A,B,\Dl
}
$$
\\ \\
$$
\infer[L\imp]{\hH~|~A\imp B,\Gm\seq\Dl}{
 \hH~|~A\imp B,\Gm\seq A,\Dl &
 \hH~|~B,\Gm\seq\Dl
}
$$
&
$$
\infer[R\imp]{\hH~|~\Gm\seq\Dl,A\imp B}{
 \hH~|~A,\Gm\seq B
}
$$
\\ \\
$$
\infer[EW]{\hH~|~\Gm\seq\Dl}{
\hH
}
$$
&
$$
\infer[EC]{\hH~|~\Gm\seq\Dl}{
 \hH~|~\Gm\seq\Dl~|~\Gm\seq\Dl
}
$$
\\ \\
where $\hH$ is non-empty in \rn{EW}
&
\end{tabular}
\caption{The hypersequent calculus \sn{HG3ipm}.  $P$ is atomic.}
\label{figure:H3ip}
\end{figure}

\begin{samepage}
\begin{proposition}[\cite{Roth2010-PhD}] The standard internal weakening and contraction rules
$$
\infer[LW]{\hH~|~A,\Gm\seq\Dl}{
 \hH~|~\Gm\seq\Dl
} \quad
\infer[RW]{\hH~|~\Gm\seq\Dl,A}{
 \hH~|~\Gm\seq\Dl
}
\quad
\infer[LC]{\hH~|~A,\Gm\seq\Dl}{
 \hH~|~A,A,\Gm\seq\Dl
} \quad
\infer[RC]{\hH~|~\Gm\seq\Dl,A}{
 \hH~|~\Gm\seq\Dl,A,A
}
$$
are admissible in \sn{HG3ipm}.
\end{proposition}
\end{samepage}

\begin{remark}[Notation]
As with labelled proofs in Remark~\ref{rem:short_struc},
multiple instances of \rn{EW}, \rn{LW} or \rn{RW} will be
shown simply as \rn{W}. Similarly, multiple instances of \rn{EC},
\rn{LC} or \rn{RC} will be shown simply as \rn{C}.
\end{remark}

We note that the many hypersequent calculi treat the components as
corresponding to points in the Kripke semantics of a logic, e.g.
\cite{Avro1991,Avro1996} or \cite{CiabFerr2001-JLC}, and use this as a
motivation for translating frame axioms into hypersequent rules. We
generalise this here by using monotonicity to encode relations between
points as subset relations between components in the following
procedure:


\begin{definition}[Translation of Geometric Axioms to Hypersequent Rules]
Structural hypersequent rules are obtained from geometric frame axioms by
the following method:\label{def:frame_axioms_to_hs}
\begin{enumerate}
\item Translate the frame axiom into a geometric rule
using the procedure from Definition~\ref{def:deriving_geo_rules}, and
 expand the sets of relations in the conclusion and premisses
  to their transitive closures.
\item Create a \defterm{base schematic hypersequent} by associating each
  principal label $x$ from the conclusion with a component containing a unique pair of multiset variables $\Gm_x$ and $\Dl_x$, e.g. for a rule with principal labels $x,y$, the base schematic hypersequent is $\hH~|~\Gm_x\seq\Dl_x~|~\Gm_y\seq\Dl_y$.
\item Take the base schematic hypersequent: for each relation $x\rel
  y$ in the conclusion, add $\Gm_x$ to the antecedent of the component
  associated with $y$, and add $\Dl_y$ to the succedent of the
  component associated with $x$.
(When there is a symmetric relation between components, they can be
merged into a single component.)
 In the previous
  example, $x\rel y$ would be
  $\hH~|~\Gm_x\seq\Dl_x,\Dl_y~|~\Gm_x,\Gm_y\seq\Dl_y$.\label{item:conc}

\item Using the result of step~\ref{item:conc}, repeat the same process for each premiss. For fresh labels, add new components, but do not add new variables.
\item Remove multiple occurrences of the same variable in the
  antecedent or succedent, as well as duplicate schematic components,
  e.g. $\Gm_x,\Gm_x,\Gm_y\seq\Gm_y$ can be changed to $\Gm_x,\Gm_y\seq\Gm_y$.
\end{enumerate}
\end{definition}

\begin{remark}\label{rem:semantic_mismatch}
We note that the treatment of the components of a hypersequent
as corresponding to points in the Kripke semantics of a logic appears
to be at odds with Definition~\ref{def:hs_semantics}.
That is equivalent to
$\exists x\in W.\mf{M},x\force A$ implies $\mf{M}\vDash A$,
which corresponds to the frame semantics for \sn{Class}.
However, by Lemma~\ref{lem:rooted}, we can assume that a model has a
distinguished point to which the distinguished component corresponds.
From monotonicity, that component is true in all points of the model.
Hence the hypersequent is true by addition.
\end{remark}

\begin{lemma} The method from Definition~\ref{def:frame_axioms_to_hs} yields
sound rules for the corresponding logics in \IntGeo.
\begin{proof} Note that Step~\ref{item:conc} constructs components that
satisfy the monotonicity property w.r.t. subsets of formulae in
the corresponding points in a model in accordance with the frame axiom.
\end{proof}
\end{lemma}

\begin{lemma} The method from Definition~\ref{def:frame_axioms_to_hs}
  yields rules which admit internal weakening and contraction in a
  \sn{HG3ipm}-like calculus.
  \begin{proof} Note that the rules have the \defterm{subformula
      property}, i.e. each multiset variable in the premisses occurs
    in the conclusion. Thus instances of internal weakening can be
    permuted to lower derivation depths.  Note also that the
    antecedents (and succedents) of components in the conclusions are
    subsets of the antecedents (and succedents) of corresponding
    components in each premiss, and that the rules are context
    sharing.  Thus instances of internal contraction can be permuted
    to lower derivation depths.
\end{proof}
\end{lemma}

Applying the method from Definition~\ref{def:frame_axioms_to_hs} to
the axioms in Table~\ref{tbl:char} yields the rules in
Fig.~\ref{fig:char_hs}.  Note that the \rn{dir}, \rn{lin} and
\rn{sym} rules are interderivable with the rules \rn{LQ},
\rn{Com} and \rn{S} from the literature,
e.g. \cite{Avro1996}. \sn{HG3ipm^*} is the system \sn{HG3ipm} augmented
by these rules.

\begin{figure}
$$
\infer[dir]{\hH | \Gm_1\seq\Dl_1 | \Gm_2\seq\Dl_2}{
 \hH | \Gm_1\seq\Dl_1 | \Gm_2\seq\Dl_2 | \Gm_1,\Gm_2\seq
}
\quad
\infer[lin]{\hH | \Gm_1\seq\Dl_1 | \Gm_2\seq\Dl_2}{
 \hH | \Gm_1\seq\Dl_1,\Dl_2 | \Gm_1,\Gm_2\seq\Dl_2 &
 \hH | \Gm_1,\Gm_2\seq\Dl_1 | \Gm_2\seq\Dl_1,\Dl_2
}
$$

$$
\infer[bd_2]{\hH~|~\Gm_1\seq\Dl_1,\Dl_2,\Dl_3~|~\Gm_1,\Gm_2\seq\Dl_2,\Dl_3~|~\Gm_1,\Gm_2,\Gm_3\seq\Dl_3}{
 \hH~|~\Gm_1,\Gm_2\seq\Dl_1,\Dl_2,\Dl_3~|~\Gm_1,\Gm_2,\Gm_3\seq\Dl_3 &
 \hH~|~\Gm_1\seq\Dl_1,\Dl_2,\Dl_3~|~\Gm_1,\Gm_2\Gm_3,\seq\Dl_2,\Dl_3
}
$$

$$
\infer[sym]{\hH~|~\Gm_1\seq\Dl_1,\Dl_2~|~\Gm_1,\Gm_2\seq\Dl_2}{
 \hH~|~\Gm_1,\Gm_2\seq\Dl_1,\Dl_2
}
$$
\caption{Hypersequent rules of well-known logics, obtained from the rules in Fig.~\ref{fig:char_rules}.}
\label{fig:char_hs}
\end{figure}

In \cite{CiabStrasTeru2009-CSL} it is shown that
for hypersequent calculi augmented with structural rules to admit the \rn{cut}
rule,
those structural rules must have \defterm{linear conclusions}---that is, the
multiset variables must not occur more than once in the conclusions
of those rules.
The method in Definition~\ref{def:frame_axioms_to_hs} yields rules
which meet this requirement, in cases where the original frame axiom
is of the form $\forall\ox.(\top\imp B)$, e.g. \rn{dir} and \rn{lin}.
For rules that do not have linear conclusions such as \rn{bd_2},
it would seem that \rn{cut} is not necessarily admissible.
However, we show below, using simply labelled calculi as an intermediate proxy
in Corollary~\ref{cor:LG3ipm*_cut}, that these rules
do admit \rn{cut} in $\sn{HG3ipm}+\rn{lin}$.

\subsection{Simply Labelled Calculi}

Simply labelled calculi such as those given in
\cite{Kang1957} and \cite{Masl1969} are
(syntactically) labelled calculi without relational formulae, but with
a similar semantics to hypersequents (see
Definition~\ref{def:sls_semantics}).  They can be treated as an alternative
notation for hypersequents, where formulae are annotated with a name
for the component in which they occur.  Translation between the two
formalisms is straightforward, and will be omitted for brevity.
The only issues with translation are in regards to a notion similar to
$\alpha$-equivalence on labels (which is addressed in
Definition~\ref{def:subset_modulo_perm}), and hypersequents with an
\defterm{empty component},
i.e. hypersequents of the form $\hH~|~\seq$. Since the empty component
is never true in any interpretation, the latter issue can be safely ignored
for logics in \IntGeo.

A simply labelled calculus \sn{LG3ipm^*} \cite{Roth2010-PhD} is given
in Fig.~\ref{figure:LG3ipm} as a translation from \sn{HG3ipm^*}.

\begin{figure}[h]
\centering
\begin{tabular}{c c}
$$
\infer[Ax]{\lb{x}{P},\nlb{\Gm}\seq\nlb{\Dl},\lb{x}{P}}{
}
$$
&
$$
\infer[L\bot]{\lb{x}{\bot},\nlb{\Gm}\seq\nlb{\Dl}}{
}
$$
\\ \\
\ifshort
\else
$$
\infer[L\land]{\nlb{\Gm},\lb{x}{A\land B}\seq\nlb{\Dl}}{
 \nlb{\Gm},\lb{x}{A},\lb{x}{B}\seq\nlb{\Dl}
}$$
&
$$
\infer[R\land]{\nlb{\Gm}\seq \lb{x}{A\land B},\nlb{\Dl}}{
 \nlb{\Gm}\seq \lb{x}{A},\nlb{\Dl} &
 \nlb{\Gm}\seq \lb{x}{B},\nlb{\Dl}
}
$$
\\ \\
\fi
$$
\infer[L\lor]{\nlb{\Gm},\lb{x}{A\lor B}\seq\nlb{\Dl}}{
 \nlb{\Gm},\lb{x}{A}\seq\nlb{\Dl} &
 \nlb{\Gm},\lb{x}{B}\seq\nlb{\Dl}
}
$$
&
$$
\infer[R\lor]{\nlb{\Gm}\seq \lb{x}{A\lor B},\nlb{\Dl}}{
 \nlb{\Gm}\seq \lb{x}{A},\lb{x}{B},\nlb{\Dl}
}
$$
\\ \\
$$
\infer[L\imp]{\lb{x}{(A\imp B)},\nlb{\Gm}\seq\nlb{\Dl}}{
 \lb{x}{(A\imp B)},\nlb{\Gm}\seq\nlb{\Dl},\lb{x}{A} &
 \lb{x}{B},\nlb{\Gm}\seq\nlb{\Dl}
}
$$
&
$$
\infer[R\imp]{\nlb{\Gm}'\seq\nlb{\Dl}',\lb{x}{\Phi},\lb{x}{(A\imp B)}}{
 \lb{x}{A},\nlb{\Gm}'\seq\nlb{\Dl}',\lb{x}{B}
}
$$
\end{tabular}
$$
\infer[dir]{\lb{x}{\Pi},\lb{y}{\Sg},\nlb{\Gm}'\seq\nlb{\Dl}}{
 \lb{zx}{\Pi},\lb{zy}{\Sg},\nlb{\Gm}'\seq\nlb{\Dl}
}
\quad
\infer[lin]{\lb{x}{\Pi},\lb{y}{\Sg},\nlb{\Gm}'\seq\nlb{\Dl}',\lb{x}{\Phi},\lb{y}{\Psi}}{
 \lb{xy}{\Pi},\lb{y}{\Sg},\nlb{\Gm}'\seq\nlb{\Dl}',\lb{x}{\Phi},\lb{xy}{\Psi} &
 \lb{x}{\Pi},\lb{xy}{\Sg}\nlb{\Gm}'\seq\nlb{\Dl}',\lb{xy}{\Phi},\lb{y}{\Psi}
}
$$

$$
\infer[bd_2]{\lb{xyz}{\Pi_1},\lb{yz}{\Pi_2},\lb{z}{\Pi_3},\nlb{\Gm}'\seq\nlb{\Dl}',\lb{x}{\Phi_1},\lb{xy}{\Phi_2},\lb{xyz}{\Phi_3}}{
 \lb{xz}{\Pi_1},\lb{xz}{\Pi_2},\lb{z}{\Pi_3},\nlb{\Gm}'\seq\nlb{\Dl}',\lb{x}{\Phi_1},\lb{x}{\Phi_2},\lb{xz}{\Phi_3} &
 \lb{xz}{\Pi_1},\lb{z}{\Pi_2},\lb{z}{\Pi_3},\nlb{\Gm}'\seq\nlb{\Dl}',\lb{x}{\Phi_1},\lb{xz}{\Phi_2},\lb{xz}{\Phi_3}
}
$$

$$
\infer[sym]{\lb{xy}{\Pi},\lb{y}{\Sg},\nlb{\Gm}'\seq\nlb{\Dl}',\lb{x}{\Phi},\lb{xy}{\Psi}}{
 \lb{x}{\Pi},\lb{x}{\Sg},\nlb{\Gm}'\seq\nlb{\Dl}',\lb{x}{\Phi},\lb{x}{\Psi}
}
$$
\caption{The calculus \sn{LG3ipm^*}. $x\#\nlb{\Dl}'$ in \rn{R\imp} and $x,y\#\nlb{\Gm}',\nlb{\Dl}'$ in the structural rules.}
\label{figure:LG3ipm}
\end{figure}

\begin{definition}\label{def:sls_semantics}

Let $\nlb{\Gm}\slice x =_{def} \{~ \lb{x}{A} ~|~ \lb{x}{A}\in\nlb{\Gm} ~\}$.
Let $\mf{M} = \struct{W,R,D}$ be a Kripke model for
a logic in \IntGeo. Then $\mf{M}\vDash\nlb{\Gm}\seq\nlb{\Dl}$
iff there exists a label $x\in\mt{lab}(\nlb{\Gm},\nlb{\Dl})$
such that $\mf{M}\vDash\nlb{\Gm}\slice x\seq\nlb{\Dl}\slice x$, i.e.
either $\mf{M}\nvDash\mland\nlb{\Gm}\slice x$ or
$\mf{M}\vDash\mlor\nlb{\Dl}\slice x$.
\end{definition}

\begin{definition}[Subset Modulo Permutation]
Let $\nlb{\Gm}\eqmp\nlb{\Dl}$ mean that two multisets of labelled formulae
are identical, modulo permutation of labels. Then
$\nlb{\Gm}\subseteqmp\nlb{\Dl}$ iff there exists $\nlb{\Gm}'$ such that
$\nlb{\Gm}'\eqmp\nlb{\Gm}$ and $\nlb{\Gm}'\subseteq\nlb{\Dl}$.
This notion is extended naturally for sequents.
\label{def:subset_modulo_perm}\end{definition}

\begin{proposition}[Label Substitution] Let $\nlb{\Gm}\seq\nlb{\Dl}$ be
a simply labelled sequent, and $x,y$ be labels.
If $\sn{LG3ipm^*}\vdash\nlb{\Gm}\seq\nlb{\Dl}$, then
$\sn{LG3ipm^*}\vdash[y/x]\nlb{\Gm}\seq[y/x]\nlb{\Dl}$.
\ifshort
\else
\begin{proof} Straightforward.
\end{proof}
\fi
\end{proposition}

\begin{proposition} Weakening and contraction
are admissible in \sn{LG3ipm^*}.
\ifshort
\else
\begin{proof} Straightforward.
\end{proof}
\fi
\end{proposition}

\section{Translation of Labelled Proofs to Simply Labelled Proofs}

\label{sec:translation}
Labelled sequents are more expressive than hypersequents. It is not
obvious what hypersequent an arbitrary labelled sequent with
relational formulae, e.g. $x\rel y;\lb{x}{A}\seq\lb{y}{A}$,
corresponds to.
We apply the idea for translating frame
axioms into hypersequent rules (Definition~\ref{def:frame_axioms_to_hs})
by using monotonicity to encode
relational formulae as subset relations between the components---this
seems to be an obvious choice.

A translation from labelled sequents to simply labelled sequents is
given below.  (The translation from simply labelled sequents to
hypersequents is straightforward, and is omitted for brevity.)

\begin{definition}[Transitive Unfolding]\label{def:trans_unfold}
Let $\Sg^+$ be the transitive closure of $\Sg$, so that
\begin{align*}
 \overrightarrow{\mt{lab}}_x(\Sg^+) &=_{def} \{~y~|~x\rel y\in\Sg^+\}
 &
 \overleftarrow{\mt{lab}}_y(\Sg^+)  &=_{def} \{~x~|~x\rel y\in\Sg^+\}
\end{align*}
Let $\overrightarrow{\Sg}$ be the \emph{list} of labels
constructed from the multiset of relational formulae $\Sg$ \cite{Szpi1930-FM},
and let $\overleftarrow{\Sg}$ be the reversed list from $\overrightarrow{\Sg}$.
Then we define the functions
\begin{align*}
\overrightarrow{\mt{TU}}\ \overrightarrow{\Sg}\ \nlb{\Gm} &=_{def}
 \begin{cases}
  \overrightarrow{\mt{TU}}\ \overrightarrow{\Sg'}\ \nlb{\Gm}\cup\bigcup (\lambda y.[y/x]\nlb{\Gm}\slice x)\ostar\overrightarrow{\mt{lab}}_x(\Sg^+) &
   \text{where $\overrightarrow{\Sg} = x::\overrightarrow{\Sg'}$} \\
  \nlb{\Gm} & \text{otherwise} \\
 \end{cases}
\\
\overleftarrow{\mt{TU}}\ \overleftarrow{\Sg}\ \nlb{\Dl} &=_{def}
 \begin{cases}
  \overleftarrow{\mt{TU}}\ \overleftarrow{\Sg'}\ \nlb{\Dl}\cup\bigcup (\lambda x.[x/y]\nlb{\Dl}\slice y)\ostar\overleftarrow{\mt{lab}}_y(\Sg^+) &
   \text{where $\overleftarrow{\Sg} = x::\overleftarrow{\Sg'}$} \\
  \nlb{\Dl} & \text{otherwise} \\
 \end{cases}
\end{align*}
where $x::\overrightarrow{\Sg'}$ and  $x::\overleftarrow{\Sg'}$ denote lists
of labels, with $x$ as the head, and
$\ostar$ is an alternative for $\mt{map}$, where
$f\ostar \Gm =_{def} \mt{map}\ f\ \Gm$
for a function $f$ and a list or multiset $\Gm$.
Then
$(\Sg;\nlb{\Gm}\seq\nlb{\Dl})^{\bullet} =_{def}
 (\overrightarrow{\mt{TU}}\ \overrightarrow{\Sg}\ \nlb{\Gm})\seq(\overleftarrow{\mt{TU}}\ \overleftarrow{\Sg}\ \nlb{\Dl})$.
\end{definition}

Note that there is no 1-1 relation between a labelled sequent and its
transitive unfolding, e.g. $(x\rel y,x\rel z;\lb{x}{A}\seq\lb{y}{B})^{\bullet}
= (x\rel y,y\rel z;\lb{x}{A}\seq\lb{y}{B})^{\bullet}$.
Furthermore, despite encoding relations between labels as subset
relations between components, there are no rules in the corresponding
simply labelled sequent (or hypersequent) calculus to preserve this
relation.  For example, take a labelled sequent
that is derivable in \sn{G3I} without any extension rules,
$x\rel y;\lb{x}{(A\lor B)},\lb{x}{(B\imp C)}\seq\lb{x}{A},\lb{y}{C}$.
It's transitive unfolding,
$\lb{xy}{(A\lor B)},\lb{xy}{(B\imp C)}\seq\lb{x}{A},\lb{xy}{C}$,
cannot be derived in \sn{LG3ipm}.
The occurrences of $A\lor B$ in two different slices must be analysed
in parallel using a rule such as Proposition~\ref{lem:l_lor_bullet}
below, which requires linearity.
This is unsurprising, as the slices (or components) correspond to
chains through points in a model, rather than points in a model.

We now show that proofs in a labelled calculus based on \sn{G3I^*} can
be translated into proofs in a simply labelled calculus based on
\sn{LG3ipm^*} for logics in \IntGeo~ augmented with the \rn{lin} rule
(which corresponds to logics based on \sn{GD}).
(We use the notation \rn{\rho_{\iota}} to indicate a ``trivially
invertible'' form of the rule \rn{\rho} with the principal formula in
all premisses, e.g. \rn{L\imp_{\iota}}. Note that the rules are
interderivable using weakening and contraction.)
The translation of proofs from the
simply labelled to hypersequent calculus \sn{HG3ipm} is straightforward, and
is omitted for brevity.

\ifshort
\else
\begin{proposition} The rule\label{prop:LG3ipm_R_bot}
$$
\infer=[R\bot]{\nlb{\Gm}\seq\nlb{\Dl}}{
 \nlb{\Gm}\seq\nlb{\Dl},\lb{x}{\bot}
}
$$
is admissible in \sn{LG3ipm^*}.
\begin{proof} By induction on the derivation depth.
\end{proof}
\end{proposition}
\fi

\begin{proposition}\label{lem:l_lor_bullet} The rule
$$
\infer[L\lor_{\bullet}]{\lb{x}{(A\lor B)},\lb{y}{(A\lor B)},\nlb{\Gm}\seq~\nlb{\Dl}}{
 \lb{x}{A},\lb{y}{A},\nlb{\Gm}\seq~\nlb{\Dl} &
 \lb{x}{B},\lb{y}{B},\nlb{\Gm}\seq~\nlb{\Dl}
}
$$
is derivable in $\sn{LG3ipm^*}+\rn{lin}$.
\ifshort
\else
\begin{proof} 
We use $\lb{x}{\Gm}_{12}$ as shorthand for $\lb{x}{\Gm}_1,\lb{x}{\Gm}_2$ below:
\begin{equation}
\vcenter{
 \infer[lin]{\lb{x}{A},\lb{y}{B},\lb{x}{\Gm}_1,\lb{y}{\Gm}_2,\nlb{\Gm}'\seq~\nlb{\Dl}',\lb{x}{\Dl}_1,\lb{y}{\Dl}_2}{
  \infer[W]{\lb{x}{A},\lb{y}{A},\lb{y}{B},\lb{x}{\Gm}_1,\lb{y}{\Gm}_{12},\nlb{\Gm}'\seq~\nlb{\Dl}',\lb{x}{\Dl}_{12},\lb{y}{\Dl}_2}{
   \lb{x}{A},\lb{y}{A},\lb{x}{\Gm}_1,\lb{y}{\Gm}_2,\nlb{\Gm}'\seq~\nlb{\Dl}',\lb{x}{\Dl}_1,\lb{y}{\Dl}_2
  }
  &
  \infer[W]{\lb{x}{A},\lb{x}{B},\lb{y}{B},\lb{x}{\Gm}_{12},\lb{y}{\Gm}_2,\nlb{\Gm}'\seq~\nlb{\Dl}',\lb{x}{\Dl}_1,\lb{y}{\Dl}_{12}}{
   \lb{x}{B},\lb{y}{B},\lb{x}{\Gm}_1,\lb{y}{\Gm}_2,\nlb{\Gm}'\seq~\nlb{\Dl}',\lb{x}{\Dl}_1,\lb{y}{\Dl}_2
  }
 }
}
\label{eq:L_lor_bullet_com_e1}\end{equation}

\begin{equation}
\vcenter{
 \infer[lin]{\lb{x}{B},\lb{y}{A},\lb{x}{\Gm}_1,\lb{y}{\Gm}_2,\nlb{\Gm}'\seq~\nlb{\Dl}',\lb{x}{\Dl}_1,\lb{y}{\Dl}_2}{
  \infer[W]{\lb{x}{B},\lb{y}{B},\lb{y}{A},\lb{x}{\Gm}_1,\lb{y}{\Gm}_{12},\nlb{\Gm}'\seq~\nlb{\Dl}',\lb{x}{\Dl}_{12},\lb{y}{\Dl}_2}{
   \lb{x}{B},\lb{y}{B},\lb{x}{\Gm}_1,\lb{y}{\Gm}_2,\nlb{\Gm}'\seq~\nlb{\Dl}',\lb{x}{\Dl}_1,\lb{y}{\Dl}_2
  }
  &
  \infer[W]{\lb{x}{B},\lb{x}{A},\lb{y}{A},\lb{x}{\Gm}_{12},\lb{y}{\Gm}_2,\nlb{\Gm}'\seq~\nlb{\Dl}',\lb{x}{\Dl}_1,\lb{y}{\Dl}_{12}}{
   \lb{x}{A},\lb{y}{A},\lb{x}{\Gm}_1,\lb{y}{\Gm}_2,\nlb{\Gm}'\seq~\nlb{\Dl}',\lb{x}{\Dl}_1,\lb{y}{\Dl}_2
  }
 }
}
\label{eq:L_lor_bullet_com_e2}\end{equation}
where $x,y\#\nlb{\Gm}',\nlb{\Dl}'$ in \eqref{eq:L_lor_bullet_com_e1}
and \eqref{eq:L_lor_bullet_com_e2}.

$$
\infer[L\lor]{\lb{x}{(A\lor B)},\lb{y}{(A\lor B)},\nlb{\Gm}\seq~\nlb{\Dl}}{
 \infer[L\lor]{\lb{x}{A},\lb{y}{(A\lor B)},\nlb{\Gm}\seq~\nlb{\Dl}}{
   \lb{x}{A},\lb{y}{A},\nlb{\Gm}\seq~\nlb{\Dl} &
   \infer*[\eqref{eq:L_lor_bullet_com_e1}]{\lb{x}{A},\lb{y}{B},\nlb{\Gm}\seq~\nlb{\Dl}}{}
 } &
 \infer[L\lor]{\lb{x}{B},\lb{y}{(A\lor B)},\nlb{\Gm}\seq~\nlb{\Dl}}{
   \infer*[\eqref{eq:L_lor_bullet_com_e2}]{\lb{x}{B},\lb{y}{A},\nlb{\Gm}\seq~\nlb{\Dl}}{
   } &
   \lb{x}{B},\lb{y}{B},\nlb{\Gm}\seq~\nlb{\Dl}
 }
}
$$
\end{proof}
\fi
\end{proposition}


\begin{proposition}[Monotonicity Rules] The rules\label{prop:subseteqmp}
$$
\infer=[L\subseteqmp]{\nlb{\Gm},\lb{x}{A}\seq\nlb{\Dl}}{
 \nlb{\Gm},\lb{xy}{A}\seq\nlb{\Dl}
}
\qquad
\infer=[R\subseteqmp]{\nlb{\Gm}\seq\nlb{\Dl},\lb{y}{A}}{
 \nlb{\Gm}\seq\nlb{\Dl},\lb{xy}{A}
}
$$
where $\nlb{\Gm}\slice x\subseteqmp\nlb{\Gm}\slice y$
and $\nlb{\Dl}\slice y\subseteqmp\nlb{\Dl}\slice x$,
are derivable in \sn{LG3ipm^*}+\rn{lin}.

\begin{proof} Let
$\nlb{\Gm} = \lb{xy}{\Pi},\lb{y}{\Sg},\nlb{\Gm'}$ and
$\nlb{\Dl} = \nlb{\Dl'},\lb{x}{\Phi},\lb{xy}{\Psi}$,
where $x,y\notin\nlb{\Gm'},\nlb{\Dl'}$. Then
$$
\infer[def]{\nlb{\Gm},\lb{xy}{A}\seq\nlb{\Dl}}{
 \infer[Com]{\lb{x}{A},\lb{xy}{\Pi},\lb{y}{\Sg},\nlb{\Gm'}\seq\nlb{\Dl'},\lb{x}{\Phi},\lb{xy}{\Psi}}{
  \infer[W]{\lb{xy}{A},\lb{xyy}{\Pi},\lb{y}{\Sg},\nlb{\Gm'}\seq\nlb{\Dl'}\lb{x}{\Phi},\lb{xxy}{\Psi}}{
   \infer[def]{\lb{xy}{A},\lb{xy}{\Pi},\lb{y}{\Sg},\nlb{\Gm'}\seq\nlb{\Dl'}\lb{x}{\Phi},\lb{xy}{\Psi}}{
    \nlb{\Gm},\lb{xy}{A}\seq\nlb{\Dl}
   }
  }
  &
  \infer[W]{\lb{x}{A},\lb{xxy}{\Pi},\lb{xy}{\Sg},\nlb{\Gm'}\seq\nlb{\Dl'},\lb{xy}{\Phi},\lb{xyy}{\Psi}}{
    \infer[C]{\lb{x}{A},\lb{x}{\Pi},\lb{x}{\Sg},\nlb{\Gm'}\seq\nlb{\Dl'},\lb{x}{\Phi},\lb{x}{\Psi}}{
      \infer[W]{\lb{xy}{A},\lb{xy}{\Pi},\lb{xy}{\Sg},\nlb{\Gm'}\seq\nlb{\Dl'},\lb{xy}{\Phi},\lb{xxyy}{\Psi}}{
   \infer[def]{\lb{xy}{A},\lb{xy}{\Pi},\lb{y}{\Sg},\nlb{\Gm'}\seq\nlb{\Dl'}\lb{x}{\Phi},\lb{xy}{\Psi}}{
    \nlb{\Gm},\lb{xy}{A}\seq\nlb{\Dl}
   }
      }
     }
   }
 }
}
$$
The inverted form of the rule is derivable using weakening.
The derivation of the rule \rn{R\subseteqmp} is similar.

\end{proof}
\end{proposition}

\ifshort
\else

\begin{proposition} The rule
$$
\infer[R\land_{\bullet}]{\nlb{\Gm}\seq\nlb{\Dl},\lb{x}{(A\land B)},\lb{y}{(A\land B)}}{
 \nlb{\Gm}\seq\nlb{\Dl},\lb{x}{A},\lb{y}{A} &
 \nlb{\Gm}\seq\nlb{\Dl},\lb{x}{B},\lb{y}{B}
}
$$
where $\nlb{\Gm}\slice x\subseteqmp\nlb{\Gm}\slice y$,
is admissible in $\sn{LG3ipm^*}+\rn{lin}$.
\begin{proof} Similar to Proposition~\ref{lem:l_lor_bullet}.
\end{proof}
\end{proposition}

\begin{corollary}\label{cor:r_and_bullet} The rule
$$
\infer[R\land_{\bullet}^*]{\nlb{\Gm}\seq\nlb{\Dl},\lb{x_1}{(A\land B)},\ldots,\lb{x_n}{(A\land B)}}{
 \nlb{\Gm}\seq\nlb{\Dl},\lb{x_1}{A},\ldots,\lb{x_n}{A} &
 \nlb{\Gm}\seq\nlb{\Dl},\lb{x_1}{B},\ldots,\lb{x_n}{B}
}
$$
where $\nlb{\Gm}\slice x_i\subseteqmp\nlb{\Gm}\slice x_n$ (for $1\le i\le n$),
is admissible in $\sn{LG3ipm^*}+\rn{lin}$.
\begin{proof} Straightforward, using \rn{R\subseteqmp} and \rn{RW}.
\end{proof}
\end{corollary}
\fi

\begin{proposition} The rule\label{lem:LG3ipm_imp_suc_elim}
$$
\infer[RC\imp]{\nlb{\Gm}\seq\lb{x}{(A\imp B)},\nlb{\Dl}}{
 \nlb{\Gm}\seq\lb{x}{B},\lb{x}{(A\imp B)},\nlb{\Dl}
}
$$
is admissible in  $\sn{LG3ipm^*}+\rn{lin}$.
\begin{proof} By induction on the derivation depth.
\end{proof}
\end{proposition}

\begin{theorem}
\label{thm:RLS_to_SLS}
Let $\Sg;\nlb{\Gm}\seq\nlb{\Dl}$ be a labelled sequent.
If $\sn{G3I^*} \vdash\Sg;\nlb{\Gm}\seq\nlb{\Dl}$,
then $\sn{LG3ipm^*}+\rn{lin} \vdash(\Sg;\nlb{\Gm}\seq\nlb{\Dl})^{\bullet}$.
\begin{proof} By induction on the derivation depth.
The proof is given in Appendix~\ref{sec:proof_thm1}.
\end{proof}
\end{theorem}

\begin{example} Take the following proof in \sn{G3I} (using
  context-splitting rules for brevity):
$$
\infer[L\lor]{x\rel y;\lb{x}{(A\lor B)},\lb{x}{(B\imp C)}\seq\lb{x}{A},\lb{y}{C}}{
 \infer[refl]{\lb{x}{A}\seq\lb{x}{A}}{
  x\rel x;\lb{x}{A}\seq\lb{x}{A}
 }
 &
 \infer[\RLS{L\imp}]{x\rel y;\lb{x}{B},\lb{x}{(B\imp C)}\seq\lb{y}{C}}{
  x\rel y;\lb{x}{B},\lb{x}{(B\imp C)}\seq\lb{y}{B}
  &
  \infer[refl]{\lb{y}{C}\seq\lb{y}{C}}{
   y\rel y;\lb{y}{C}\seq\lb{y}{C}
  }
 }
}
$$
From Theorem~\ref{thm:RLS_to_SLS}, we can construct a proof of $\TU{\big(x\rel y;\lb{x}{(A\lor B)},\lb{x}{(B\imp C)}\seq\lb{x}{A},\lb{y}{C}\big)}$ in
\sn{LG3ipm^*}:
$$
\infer[L\lor_{\bullet}]{\lb{xy}{(A\lor B)},\lb{xy}{(B\imp C)}\seq\lb{x}{A},\lb{xy}{C}}{
 \infer[C]{\lb{xy}{A}\seq\lb{x}{A}}{
  \lb{xxy}{A}\seq\lb{x}{A}
 }
 &
 \infer[L\imp_{\iota}]{\lb{xy}{B},\lb{xy}{(B\imp C)}\seq\lb{xy}{C}}{
  \infer[R\subseteq]{\lb{xy}{B},\lb{xy}{(B\imp C)}\seq\lb{y}{B}}{
   \lb{xy}{B},\lb{xy}{(B\imp C)}\seq\lb{xy}{B}
  }
  &
  \infer[C]{\lb{y}{C}\seq\lb{xy}{C}}{
   \lb{yy}{C}\seq\lb{xy}{C}
  }
 }
}
$$
Note that the contractions are superfluous for this example.
\end{example}

\begin{theorem}
\label{thm:RLS_to_SLS_cut}
Let $\Sg;\nlb{\Gm}\seq\nlb{\Dl}$ be a labelled sequent.
If $\sn{G3I^*}+\rn{cut} \vdash\Sg;\nlb{\Gm}\seq\nlb{\Dl}$,
then $\sn{LG3ipm^*}+\rn{lin}+\rn{cut} \vdash(\Sg;\nlb{\Gm}\seq\nlb{\Dl})^{\bullet}$.
\begin{proof} By induction on the derivation depth, similar to
the proof of Theorem~\ref{thm:RLS_to_SLS}, with an additional case for
the \rn{cut} rule.
Suppose that $\Sg,\Sg';\nlb{\Gm},\nlb{\Gm}'\seq\nlb{\Dl},\nlb{\Dl}'$
is the conclusion of an instance of \rn{cut}:
$$
\infer[cut]{\Sg,\Sg';\nlb{\Gm},\nlb{\Gm}'\seq\nlb{\Dl},\nlb{\Dl}}{
 \Sg;\nlb{\Gm}\seq\nlb{\Dl},\lb{x}{A} &
 \Sg';\lb{x}{A},\nlb{\Gm}'\seq\nlb{\Dl}'
}
$$
Let
\begin{align*}
\TU{(\Sg;\nlb{\Gm}\seq\nlb{\Dl},\lb{x}{A})} &= \TU{\Gm}\seq\TU{\Dl},\lb{y_1}{A},\ldots,\lb{y_m}{A},\lb{x}{A} \\
\TU{(\Sg';\lb{x}{A},\nlb{\Gm}'\seq\nlb{\Dl}')} &= \lb{z_1}{A},\ldots,\lb{z_n}{A},\lb{x}{A},\TU{\Gm'}\seq\TU{\Dl'} \\
\TU{\Sg,\Sg';\nlb{\Gm},\nlb{\Gm}'\seq\nlb{\Dl},\nlb{\Dl}} &= \TU{\Gm},\TU{\Gm'}\seq\TU{\Dl},\TU{\Dl'}
\end{align*}
where $\TU{\Dl}\slice x \subseteqmp \TU{\Dl}\slice y_i$ (for $1\le i\le m$)
and $\TU{\Gm'}\slice x \subseteqmp \TU{\Gm'}\slice z_i$ (for $1\le i\le n$).
Then the corresponding proof in $\sn{LG3ipm^*}+\rn{lin}+\rn{cut}$ is derived:
$$
\infer[cut]{\TU{\Gm},\TU{\Gm'}\seq\TU{\Dl},\TU{\Dl'}}{
 \infer[R\subseteqmp^m]{\TU{\Gm},\TU{\Gm'}\seq\TU{\Dl},\TU{\Dl'},\lb{x}{A}}{
  \TU{\Gm}\seq\TU{\Dl},\lb{y_1}{A},\ldots,\lb{y_m}{A},\lb{x}{A}
 }
 &
 \infer[L\subseteqmp^n]{\lb{x}{A},\TU{\Gm'}\seq\TU{\Dl'}}{
  \lb{z_1}{A},\ldots,\lb{z_n}{A},\lb{x}{A},\TU{\Gm'}\seq\TU{\Dl'}
 }
}
$$
\end{proof}
\end{theorem}

\begin{corollary}\label{cor:LG3ipm*_cut}
$\sn{LG3ipm^*}+\rn{lin}$ admits \rn{cut}.
\begin{proof} Follows from Lemma.~\ref{prop:G3I_struc}
and Theorem~\ref{thm:RLS_to_SLS_cut}. We note that in cases where
cuts can be eliminated in \sn{G3I^*}, they can be eliminated in
the corresponding
translations into $\sn{LG3ipm^*}+\rn{lin}$.
\end{proof}
\end{corollary}

\section{Discussion}

We gave a method of translating relational frame axioms for logics in
\IntGeo~ into structural rules for hypersequent calculi
(Definition~\ref{def:frame_axioms_to_hs}) by
encoding the monotonicity property of those logics as subset
relations.  The resulting rules can be shown to admit weakening
and contraction in a straightforward manner.

We presented an alternative notation for hypersequents, called
simply labelled sequents, and then introduced a similar
translation method called Transitive Unfolding (Definition~\ref{def:trans_unfold})
to translate labelled proofs into simply labelled proofs
(Theorems~\ref{thm:RLS_to_SLS} and~\ref{thm:RLS_to_SLS_cut})
for logics based on \sn{GD} instead of \sn{Int}. (Recall that
$\sn{GD} = \sn{Int} + \rn{lin}$.)

The resulting simply labelled (equivalent to hypersequent) rules admit
\rn{cut} in the presence of the \rn{lin} rule.
Cut admissibility would appear surprising in light of
results in \cite{CiabGalaTeru2008-IEEESLCS}, because some of the
resulting rules have ``non-linear'' conclusions---that is, conclusions
where some metavariables occur more than once as a result of
transitive unfolding.  However, the
\rn{lin} rule allows us to derive the ``Monotonicity Rules''
(Proposition~\ref{prop:subseteqmp}) that allow us to
eliminate the duplicate metavariables, because they occur
in the correct configuration from transitive unfolding.
Instances of \rn{cut} can be permuted to the leaves of proof and eliminated
(literally by translating the corresponding labelled proof).
Thus it is not inconsistent with results in
\cite{CiabGalaTeru2008-IEEESLCS}.

That the translations require a stronger logic is not surprising,
considering that labelled sequents are more expressive than simply
labelled sequents and hypersequents.  But this is also
disappointing, since the translated proofs are not in the
original logic---for example, labelled proofs in \sn{BD_2} are
translated into simply labelled/hypersequent proofs in \sn{G_3}
(three-valued G\"odel Logic \cite{Gode1933-EMK}).

It is also noteworthy that \cite{CiabFerr2000-TABLEAUX} provides an
alternative hypersequent calculus for logics such as \sn{BD_2} by
restricting the external permutation rule, thus making the hypersequents
linear, instead of requiring a subset relationship between components
with implied linearity.

Further investigation may show in what cases, if any, the \rn{lin}
rule can be eliminated from such proofs, and may make explicit the
expressive limits of simply labelled and hypersequent calculi.


We note that this work can be adapted for similar labelled calculi,
such as the intuitionistic fragment of the system given in
\cite{PintUust2009,DBLP:journals/corr/abs-1101-5445}.
This work also can be easily adapted to
calculi for other families of logics, such as modal logics,
so long as they are {normal} logics---that is, they have
a preordered and monotonic relational semantics.

We have omitted an explicit discussion
on translating labelled calculi into hypersequent calculi, although we
believe that the method for translating geometric rules into
structural hypersequent rules can also be adapted to logical rules as
well.  It is an area for future investigation.


\section{Acknowledgements}

We are grateful to Roy Dyckhoff and Sara Negri for providing manuscripts of
earlier drafts of \cite{DyckNegr2011}, and for comments from anonymous referees.

\bibliographystyle{plain}
\bibliography{macros,logic,lattice,rothenberg}

\appendix

\section{Proof of Theorem \ref{thm:RLS_to_SLS}}\label{sec:proof_thm1}

\begin{theorem*}
Let $\Sg;\nlb{\Gm}\seq\nlb{\Dl}$ be a labelled sequent.\\
If $\sn{G3I^*} \vdash\Sg;\nlb{\Gm}\seq\nlb{\Dl}$,
then $\sn{LG3ipm^*}+\rn{lin} \vdash(\Sg;\nlb{\Gm}\seq\nlb{\Dl})^{\bullet}$.
\begin{proof} By induction on the derivation depth.

\begin{enumerate}

\ifshort
\else
\item Suppose $\Sg;\nlb{\Gm}\seq\nlb{\Dl}$ is an axiom. Then $(\Sg;\nlb{\Gm}\seq\nlb{\Dl})^{\bullet}$ is
also an axiom.



\item Suppose $\Sg;\nlb{\Gm}\seq\nlb{\Dl}$ is the conclusion of an instance of
  \rn{refl}.  We apply an instance of contraction to remove duplicate
  formulae labelled with $x$ in the antecedent and succedent.

\item Suppose $\Sg;\nlb{\Gm}\seq\nlb{\Dl}$ is the conclusion of an instance of
  \rn{trans}.  We apply an instance of contraction to remove
  duplicate formulae labelled with $z$ (unfolded from $x$) in the
  antecedent, and to remove duplicate formulae labelled with $x$
  (unfolded from $z$) in the succedent.


\item Suppose $\Sg;\nlb{\Gm}\seq\nlb{\Dl}$ is the conclusion of an instance of \rn{L\land}:
$$
\infer[L\land]{\Sg;\nlb{\Gm},\lb{x_1}{(A\land B)}\seq\nlb{\Dl}}{
 \Sg;\nlb{\Gm},\lb{x_1}{A},\lb{x_1}{B}\seq\nlb{\Dl}
}
$$
Let $(\Sg;\nlb{\Gm},\lb{x_1}{A},\lb{x_1}{B}\seq\nlb{\Dl})^{\bullet} = \nlb{\Gm}^{\bullet},\lb{x_1}{A},\lb{x_1}{B},\ldots,\lb{x_n}{A},\lb{x_n}{B}\seq\nlb{\Dl}^{\bullet}$, where $x_1\rel x_i\in\Sg^+$ ($2\le i\le n$).
The corresponding proof in \sn{LG3ipm^*} is derived using $n$ instances of
\rn{L\land}:
$$
\infer[L\land^n]{\nlb{\Gm}^{\bullet},\lb{x_1}{(A\land B)},\ldots,\lb{x_n}{(A\land B)}\seq\nlb{\Dl}^{\bullet}}{
 \nlb{\Gm}^{\bullet},\lb{x_1}{A},\lb{x_1}{B},\ldots,\lb{x_n}{A},\lb{x_n}{B}\seq\nlb{\Dl}^{\bullet}
}
$$
\label{item:rls2sls_l_land}

\item Suppose $\Sg;\nlb{\Gm}\seq\nlb{\Dl}$ is the conclusion of an instance of \rn{R\land}:
$$
\infer[R\land]{\Sg;\nlb{\Gm}\seq\lb{x_1}{(A\land B)},\nlb{\Dl}}{
 \Sg;\nlb{\Gm}\seq\lb{x_1}{A},\nlb{\Dl} &
 \Sg;\nlb{\Gm}\seq\lb{x_1}{B},\nlb{\Dl}
}
$$
where $x_2\rel x_1,\ldots,x_n\rel x_1\in\Sg^+$ for $n\ge 1$. Let
\begin{align*}
(\Sg;\nlb{\Gm}\seq\lb{x_1}{A},\nlb{\Dl})^{\bullet} &=
 \nlb{\Gm}^{\bullet}\seq\lb{x_1}{A},\ldots,\lb{x_n}{A},\nlb{\Dl}^{\bullet} \\
(\Sg;\nlb{\Gm}\seq\lb{x_1}{B},\nlb{\Dl})^{\bullet} &=
 \nlb{\Gm}^{\bullet}\seq\lb{x_1}{B},\ldots,\lb{x_n}{B},\nlb{\Dl}^{\bullet}
\end{align*}
where
  $\nlb{\Gm}\slice x_1\subseteqmp\nlb{\Gm}\slice x_i$ and
  $\nlb{\Dl}\slice x_i\subseteqmp\nlb{\Dl}\slice x_1$ for $1\le i\le n$.
The corresponding proof in \sn{LG3ipm^*} is derived using the
\rn{R\land^{\bullet}} rule from Proposition~\ref{cor:r_and_bullet}:
$$
\infer[R\land^{\bullet}]{\nlb{\Gm}^{\bullet}\seq\lb{x_1}{(A\land B)},\ldots,\lb{x_n}{(A\land B)},\nlb{\Dl}^{\bullet}}{
 \nlb{\Gm}^{\bullet}\seq\lb{x_1}{A},\ldots,\lb{x_n}{A},\nlb{\Dl}^{\bullet} &
 \nlb{\Gm}^{\bullet}\seq\lb{x_1}{B},\ldots,\lb{x_n}{B},\nlb{\Dl}^{\bullet}
}
$$
Note that $(\Sg;\nlb{\Gm}\seq\lb{x_1}{(A\land B)},\nlb{\Dl})^{\bullet} =
 \nlb{\Gm}^{\bullet}\seq\lb{x_1}{(A\land B)},\ldots,\lb{x_n}{(A\land B)},\nlb{\Dl}^{\bullet}$.

\label{item:rls2sls_r_land}

\item Suppose $\Sg;\nlb{\Gm}\seq\nlb{\Dl}$ is the conclusion of an instance of \rn{L\lor}.
The case is the dual of case~\ref{item:rls2sls_r_land} above, using the
\rn{L\lor_{\bullet}} rule from Proposition~\ref{lem:l_lor_bullet}.
\label{item:rls2sls_r_imp}

\item Suppose $\Sg;\nlb{\Gm}\seq\nlb{\Dl}$ is the conclusion of an instance of \rn{R\lor}.
The case is the dual of case~\ref{item:rls2sls_l_land} above.
\fi
\item Suppose $\Sg;\nlb{\Gm}\seq\nlb{\Dl}$ is the conclusion of an instance of \rn{\RLS{L\imp}}:
$$
\infer[\RLS{L\imp}]{x_1\rel y_1\Sg;\lb{x_1}{(A\imp B)},\nlb{\Gm}\seq\nlb{\Dl}}{
 x_1\rel y_1\Sg;\lb{x_1}{(A\imp B)},\nlb{\Gm}\seq\nlb{\Dl},\lb{y_1}{A} &
 x_1\rel y_1\Sg;\lb{x_1}{(A\imp B)},\lb{y_1}{B},\nlb{\Gm}\seq\nlb{\Dl}
}
$$
where $x_1\rel y_1,\ldots, x_m\rel y_1,y_1\rel y_2,\ldots,y_1\rel y_n$ for
$m,n\ge 1$. Let
\begin{align}
 (x_1\rel y_1\Sg;\lb{x_1}{(A\imp B)},\nlb{\Gm}\seq\nlb{\Dl},\lb{y_1}{A})^{\bullet} &=
   \overline{(A\imp B)},\TU{\nlb{\Gm}}\seq\TU{\nlb{\Dl}},\lb{x_1}{A},\ldots,\lb{x_m}{A},\lb{y_1}{A}\label{L_imp_p1_equiv}\\
 (x_1\rel y_1\Sg;\lb{x_1}{(A\imp B)},\lb{y_1}{B},\nlb{\Gm}\seq\nlb{\Dl})^{\bullet} &=
   \overline{(A\imp B)},\lb{y_1}{B},\ldots,\lb{y_n}{B},\TU{\nlb{\Gm}}\seq\TU{\nlb{\Dl}}\label{L_imp_p2_equiv}
\end{align}
where $\overline{(A\imp B)} = \lb{x_1}{(A\imp B)},\lb{y_1}{(A\imp B)},\lb{y_n}{(A\imp B)}$.
We can derive the following from \eqref{L_imp_p1_equiv}, for $1\le i\le n$:
$$
\infer[R\subseteqmp]{\overline{(A\imp B)},\TU{\nlb{\Gm}}\seq\TU{\nlb{\Dl}},\lb{y_i}{A}}{
 \infer[RW]{\overline{(A\imp B)},\TU{\nlb{\Gm}}\seq\TU{\nlb{\Dl}},\lb{y_1}{A},\lb{y_i}{A}}{
  \infer[R\subseteqmp^+]{\overline{(A\imp B)},\TU{\nlb{\Gm}}\seq\TU{\nlb{\Dl}},\lb{y_1}{A}}{
   \overline{(A\imp B)},\TU{\nlb{\Gm}}\seq\TU{\nlb{\Dl}},\lb{x_1}{A},\ldots,\lb{x_m}{A},\lb{y_1}{A}
  }
 }
}
$$
(Clearly the last two inference steps are omitted for $i=1$.)
We first derive the following:
\begin{equation}
\vcenter{
 \infer[L\imp_{\iota}]{\overline{(A\imp B)},\lb{y_2}{B},\ldots,\lb{y_n}{B},\TU{\nlb{\Gm}}\seq\TU{\nlb{\Dl}}}{
  \infer[W]{\overline{(A\imp B)},,\lb{y_2}{B},\ldots,\lb{y_n}{B},\TU{\nlb{\Gm}}\seq\TU{\nlb{\Dl}},\lb{y_1}{A}}{
   \overline{(A\imp B)},\TU{\nlb{\Gm}}\seq\TU{\nlb{\Dl}},\lb{y_1}{A}
  }
  &
  \infer*[\eqref{L_imp_p2_equiv}]{\hspace{1cm}}{}
 }
}
\label{L_imp_e1_equiv}\end{equation}
For $n\ge 2$, we apply the result of \eqref{L_imp_e1_equiv} to
$$
\infer[L\imp_{\iota}]{\overline{(A\imp B)},\lb{y_{i+1}}{B},\ldots,\lb{y_n}{B},\TU{\nlb{\Gm}}\seq\TU{\nlb{\Dl}}}{
  \infer[W]{\overline{(A\imp B)},,\lb{y_{i+1}}{B},\ldots,\lb{y_n}{B},\TU{\nlb{\Gm}}\seq\TU{\nlb{\Dl}},\lb{y_i}{A}}{
  \overline{(A\imp B)},\TU{\nlb{\Gm}}\seq\TU{\nlb{\Dl}},\lb{y_2}{A}
 }
 &
 \infer*[]{\overline{(A\imp B)},\lb{y_i}{B},\ldots,\lb{y_n}{B},\TU{\nlb{\Gm}}\seq\TU{\nlb{\Dl}}}{}
}
$$
and apply repeatedly until we have derived
$\overline{(A\imp B)},\TU{\nlb{\Gm}}\seq\TU{\nlb{\Dl}}$.

\item Suppose $\Sg;\nlb{\Gm}\seq\nlb{\Dl}$ is the conclusion of an instance of \rn{\RLS{R\imp_{\iota}}}:
$$
\infer[\rn{\RLS{R\imp_{\iota}}}]{\Sg;\lb{x_1}{\Gm},\nlb{\Gm}'\seq\nlb{\Dl}',\lb{x_1}{\Dl},\lb{x_1}{(A\imp B)}}{
  x_1\rel y,\Sg;\lb{x_1}{\Gm},\lb{y}{A},\nlb{\Gm}'\seq\nlb{\Dl}',\lb{y}{B},\lb{x_1}{\Dl},\lb{x_1}{(A\imp B)}
}
$$
where $x_2\rel x_1,\ldots,x_n\rel x_1\in\Sg^+$ for $n\ge 1$.
Let
\begin{align*}
(x_1\rel y,\Sg;&\lb{x_1}{\Gm},\lb{y}{A},\nlb{\Gm}'\seq\nlb{\Dl}',\lb{y}{B},\lb{x_1}{\Dl},\lb{x_1}{(A\imp B)})^{\bullet} = \\
 &\nlb{\Gm}^{\bullet},\lb{y}{\Gm},\lb{y}{A}\seq\nlb{\Dl}^{\bullet},\lb{y}{B},\lb{x_1}{B},\ldots,\lb{x_n}{B},\lb{x_1}{(A\imp B)},\ldots,\lb{x_n}{(A\imp B)}
\end{align*}
where $\lb{y}{\Gm}\eqmp\nlb{\Gm}^{\bullet}\slice x_1$.
The corresponding proof in \sn{LG3ipm^*} is derived:
$$
\infer[R\imp_{\iota}]{\nlb{\Gm}^{\bullet}\seq\nlb{\Dl}^{\bullet},\lb{x_1}{(A\imp B)},\ldots,\lb{x_n}{(A\imp B)}}{
 \infer[RC\imp^n]{\nlb{\Gm}^{\bullet},\lb{y}{\Gm},\lb{y}{A}\seq\nlb{\Dl}^{\bullet},\lb{y}{B},\lb{x_1}{(A\imp B)},\ldots,\lb{x_n}{(A\imp B)}}{
  \nlb{\Gm}^{\bullet},\lb{y}{\Gm},\lb{y}{A}\seq\nlb{\Dl}^{\bullet},\lb{y}{B},\lb{x_1}{B},\ldots,\lb{x_n}{B},\lb{x_1}{(A\imp B)},\ldots,\lb{x_n}{(A\imp B)}
  }
}
$$

\ifshort
\else

\item Suppose $\Sg;\nlb{\Gm}\seq\nlb{\Dl}$ is the conclusion of an
  instance of ordering rules such as \rn{dir}, \rn{lin} or \rn{sym}.
  The corresponding proof in \sn{LG3ipm^*} is derived using the simply
  labelled form of that rule, with weakening and contraction as
  appropriate. (Recall the method for deriving the hypersequent rule
  from the corresponding geometric rule.)

\fi

\end{enumerate}
\end{proof}
\end{theorem*}

\end{document}

%% file: macros.tex
\usepackage{ifthen}

\newif\ifpdf
\ifx\pdfoutput\undefined
  \pdffalse
\else
  \ifx\pdfoutput\relax
  \else
    \ifnum\pdfoutput>0
      \pdftrue
    \fi
  \fi
\fi

\newif\ifdraft 
\newif\ifshort 

\newif\iffalse 
\falsefalse

\newenvironment{notation*}{\par\noindent\bf Notation:~\rm}{}

\newcounter{mboldlevel}
\let\oldmb\mathbf
 \renewcommand{\mathbf}[1]{
   \addtocounter{mboldlevel}{1}
   \oldmb{#1}
   \addtocounter{mboldlevel}{-1}
}

\newenvironment{leftalign*}{%
\begin{equation*}\begin{array}{l l r}
}{\end{array}\end{equation*}}

\newcommand{\mb}[1]{\ensuremath{\mathbf{#1}}}
\newcommand{\mc}[1]{\ensuremath{\mathcal{#1}}}
\newcommand{\mt}[1]{\ensuremath{\mathtt{#1}}}
\newcommand{\mf}[1]{\ensuremath{\mathfrak{#1}}}

\newcommand{\ms}[1]{\lower -0.3ex\hbox{\ensuremath{\scriptstyle #1}}}

\providecommand{\url}[1]{$\langle$\texttt{#1}$\rangle$}
\providecommand{\doi}[1]{\url{doi:#1}}

\newcommand{\sn}[1]{\mb{#1}}

\newcommand{\defterm}[1]{\textbf{#1}}

\newcommand{\Gm}{\Gamma}
\newcommand{\Dl}{\Delta}
\newcommand{\Sg}{\Sigma}

\newcommand{\hs}[1]{\mc{#1}}         
\newcommand{\hG}{\hs{G}}                  
\newcommand{\hH}{\hs{H}}                  

\newcommand{\p}{\ifmmode ^{\prime} \else \(^{\prime}\) \fi}
\newcommand{\pp}{\ifmmode ^{\prime\prime}
  \else \(^{\prime\prime}\) \fi}
\newcommand{\ppp}{\ifmmode ^{\prime\prime\prime}
  \else \(^{\prime\prime\prime}\) \fi}

\newcommand{\struct}[1]{\ensuremath{\langle #1\rangle}}

\newcommand{\ostar}{\circledast}

\newcommand{\set}[1]{{#1}}

\newcommand{\imp}{\!\supset}                

\newcommand{\seq}{\!\Rightarrow\!}            

\usepackage{slashed}
\def\mland{\declareslashed{}{\wedge}{0.4}{0}{\wedge}\slashed{\wedge}\,}
\def\mlor{\declareslashed{}{\vee}{0.4}{0}{\vee}\slashed{\vee}\,}

\newcommand{\rel}{\le}                       

\newcommand{\force}{\!\Vdash\!}
\newcommand{\nforce}{\!\nVdash\!}

\newcommand{\lb}[2]{#2^{#1}}

\newcommand{\nlb}[1]{\underline{#1}}

\newcommand{\slice}{\sslash}    

\newcommand{\HIDE}[1]{}

\newcommand{\eqmp}{\approx}   

\newcommand{\subseteqmp}{\mathbin{\raisebox{0.2ex}{\ensuremath{\subset}}\raisebox{-1ex}{\ensuremath{\mkern-11mu\sim}}}}

\newcommand{\Luka}{{\L}uka\-sie\-wicz }